# О ДИСКРЕТНЫХ ГРУППАХ СИММЕТРИИ ПОДОБИЯ ЕВКЛИДОВА ПРОСТРАНСТВА

**Александр С. Прохода, к.ф.-м.н.**

Украина

Отправляясь от классических результатов Шубникова и Заморзаева, реализованы компьютерные модели фигур, позволяющие визуализировать действие дискретных подгрупп непрерывных топологических групп. Визуализация действия осуществляется, путем выполнения разбиений фигур на фундаментальные области дискретных групп симметрии подобия, определение которых было дано Заморзаевым. Особое внимание уделено моделям квазирешеток, при помощи которых построены такие разбиения фигур, что многоцветные их раскраски позволяют исследовать, в том числе и действие цветных групп симметрии подобия. Показано, что для двумерных квазикристаллов (квазирешеток) коэффициентами гомотетии являются целые алгебраические числа квадратичного расширения поля рациональных чисел. Проведено исследование некоторых конформных отображений кристаллических множеств, что позволило установить соответствие между стереоциклической проекцией и системой координат обратных радиусов.

\_\_\_\_\_\_\_\_\_\_\_\_\_\_\_\_\_\_\_\_\_

e-mail: a-prokhoda@mail.ru

# I. ПРЕОБРАЗОВАНИЕ ПОДОБИЯ

Исторически в далеком прошлом ученые, художники и архитекторы применяли в своих работах понятие о симметрии подобия. Леонардо да Винчи был одним из первых, сознательно прилагавших симметрию подобия под видом перспективы. Он же пользовался при создании разнообразных конструкций построением, которое было названо им золотым сечением [1]. Интересовался законами подобия также Гёте, видя в закрученных формах растений и животных важную для природы спиральную тенденцию [2] (см. рис. 0). Особо следует выделить работы М. Эшера, который при написании своих картин применял элементы геометрии Лобачевского (напомним, что в геометрии Лобачевского преобразование подобия является движением).

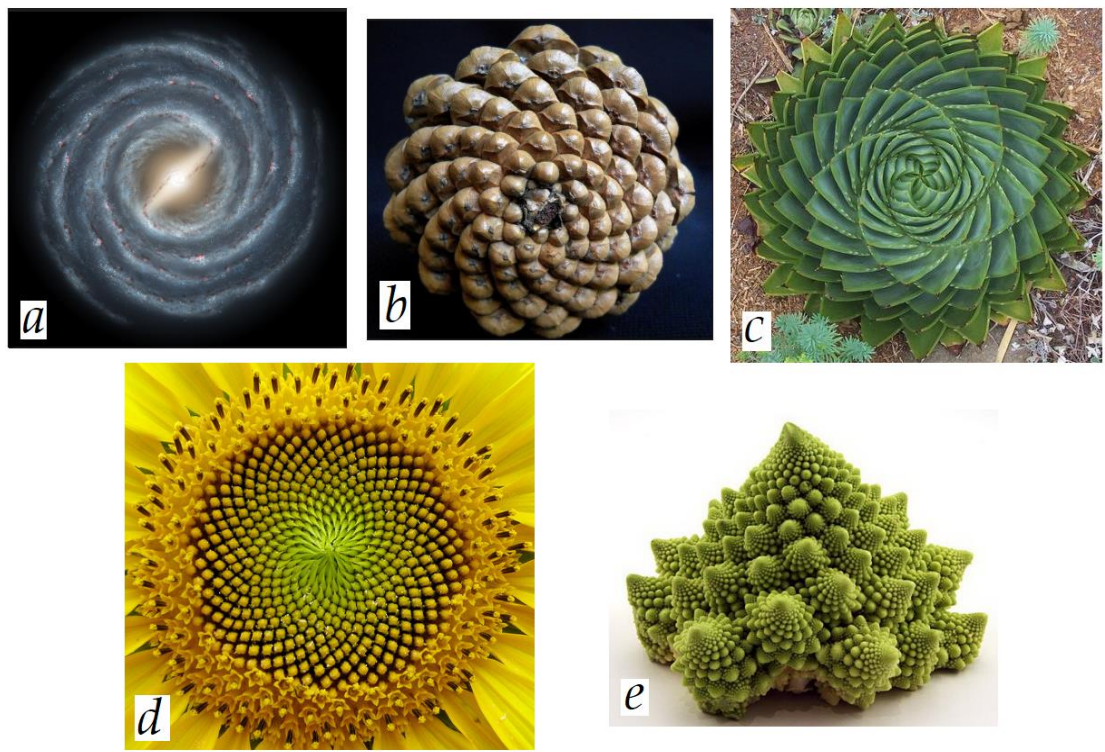

Рис. 0. Симметрия в природе (изображения взяты из Google images). *a* – спиральная галактика; *b* – сосновая шишка; *c* – разновидность алоэ; *d* – цветок подсолнечника; *e* – капуста Романеско.

Мадисон в работе [3]: ввел в рассмотрение идеализированный объект – бесконечно измельченное разбиение Пенроуза; показал, что его операциями симметрии являются операции подобия; дал положительный ответ на вопрос о возможности композиции операций подобия с различными особыми точками; установил, что преобразования апериодических кристаллов, сохраняющие ориентацию, изоморфны некоторой дискретной подгруппе группы Мёбиуса $PSL(2, C)$, т. е. могут быть реализованы как дискретные подгруппы полной группы движений пространства Лобачевского; построил стереографическую проекцию семейства локсодром на комплексную плоскость.

Некоторые конкретные случаи симметрии подобия описываются строго математически при помощи непрерывных дробей и, в частности, приводились к уже упомянутому золотому сечению. Помимо числа

$$\tau = \tau^2 - 1 = \tau^{-1} + 1 = (2\sin(18°))^{-1} = 2\cos(36°) = 2\sin(54°) = \frac{1+\sqrt{5}}{2} = 1.618034 \ ...$$

известны и другие целые алгебраические числа, именуемые металлическими сечениями. Такие числа применяются для математического моделирования двумерных и трехмерных квазирешеток (см. например работу [4]).

Множество квазикристаллов рассматриваемое, как подкласс класса кристаллов, характеризуется тем, что квазикристаллы обладают ориентационным дальним порядком. Это означает, что разбиение типа Пенроуза, например (P3), моделирующее двумерный квазикристалл обладает локальным изоморфизмом Конвея состоящим в том, что это разбиение не является периодическим, однако любой конечный фрагмент встречается в нем бесконечное число раз и обязательно появляется в круге (в случае трехмерного пространства – в шаре) достаточно большого радиуса с любым центром в этом пространстве. С другой стороны дифрактограмма квазикристалла обладает приближенной масштабной инвариантностью или скейлингом, т.е. она может быть получена однократным или многократным увеличением (уменьшением) своей собственной центральной части.

Отметим, что в случае классического кристалла его дифракционную картину определяет обратная (взаимная или дуальная) решетка. Приведем известное высказывание Бюргера [5] «Я могу с уверенностью сказать, что взаимные решетки дают нам одно из самых важных орудий для исследования дифракции рентгеновских лучей на кристаллах».

В аффинной геометрии, имеются преобразования, приводящие к масштабной инвариантности, а именно: равномерное сжатие к прямой или плоскости с коэффициентом сжатия $k$; гомотетия с центром в точке $O$ и коэффициентом гомотетии $k$. В случае скейлинга также изменяется расстояние между точками, но отношения между расстояниями сохраняются.

Приведем определение гомотетии в вещественном точечно-векторном пространстве.

**Определение 1.** Гомотетией с центром в точке $O$ и коэффициентом гомотетии $k$ называется такое отображение $K$ пространства, при котором вектор $\vec{OM}$ с началом в точке $O$ и концом в точке $M$ переводится в вектор $\vec{OM'}$ такой, что вектор $\vec{OM'} = k\vec{OM}$.

Если $K$ гомотетия с коэффициентом гомотетии $k$, то $K^{-1}$ является гомотетией с коэффициентом $k^{-1}$ ($K^{-1}$ – обратное отображение). Очевидно, что $KK^{-1} = \text{Id}$, где $\text{Id}$ – тождественное отображение. Все степени гомотетии $K$ образуют бесконечную циклическую группу (относительно операции последовательного выполнения гомотетий), которая изоморфна группе целых чисел $Z$. Наличие дальнего порядка означает наличие бесконечной группы. Например, как отмечает Кокстер в [6], изучение дифракции рентгеновских лучей показывает, что симметричная структура кристаллов является следствием симметричного расположения атомов, молекул и т.д.

Иначе говоря, имеется бесконечная группа симметрий, переводящая узор из атомов в себя (при этом считается, что этот узор распространяется на все пространство: случай идеального кристалла). В случае обычных кристаллов – это решетка: свободная абелева группа ранга равного размерности пространства. В случае квазикристаллов такую роль играет дискретная группа, являющаяся подгруппой группы преобразований подобия. Гомотетия – это частный случай преобразования подобия определяемого следующим образом.

**Определение 2.** Преобразованием подобия евклидового пространства называется такое отображение точек этого пространства, при котором сохраняется отношение расстояний между точками, это отношение называется коэффициентом подобия.

Если через $d(x, y)$ обозначить расстояние между точками на плоскости или в пространстве (в последнем случае точка $X = (x_1, x_2, …, x_n)$, $Y = (y_1, y_2, …, y_n)$ – последовательность из $n$ чисел), то при преобразовании подобия имеем, что $d(x, y) = Kd(h(x), h(y))$, $K > 0$, где $h(x), h(y)$ образы точек $x$ и $y$ при отображении $h$.

Всякое преобразование подобия является произведением (последовательным выполнением) гомотетии и движения, при этом порядок их выполнения безразличен, т.е. они коммутируют. Все преобразования подобия образуют группу, которая называется группой преобразования подобия. Эта группа является подгруппой 12-параметрической группы Aff (*V*) всех аффинных преобразований пространства *V*. Аффинные преобразования пространства характеризуются тем, что они переводят прямые линии в прямые, а преобразование подобия – это преобразование, сохраняющее форму фигуры. В классификации геометрий в соответствии с Эрлангенской программой предложенной в 1872 году Феликсом Клейном [7] аффинная группа обозначается $\mathfrak{G}_{12}$, а ее подгруппа преобразований подобия, которую, Клейн называет <u>главной группой</u>, обозначается $\mathfrak{S}_7$. Она зависит от 7 параметров. Действительно те аффинные преобразования кроме бесконечно далекой плоскости проективной группы $\mathfrak{S}_{15}$ переводят в себя и также мнимую окружность сфер и поэтому накладывает на 12 констант аффинной группы $6 - 1 = 5$ условий, так что остается как раз $12 - 5 = 7$ параметров группы $\mathfrak{S}_7$. Дискретные подгруппы групп $\mathfrak{S}_{12}$ и $\mathfrak{S}_7$ были введены в рассмотрение значительно позже.

Определение аффинной кристаллографической группы приведено в книге Артамонова и Словохотова [8]: подгруппа этой группы будет группа, которую ввел в рассмотрение Шубников [2] и назвал ее группой симметрии подобия. Причем Шубников назвал симметрию подобия – преобразованием, сохраняющим форму фигуры, а привел примеры дискретных групп преобразований подобия, в которых коэффициенты подобия являются целыми числами. Точное определение группы симметрии подобия было дано Заморзаевым в работе [9]. Приведем это определение.

**Определение 3.** *n*-мерной группой симметрии подобия называется группа симметрии подобия в *n*-мерном евклидовом пространстве обладающая следующими свойствами: 1) в ней содержится хоть одно преобразование $P$ с коэффициентом $k \neq 1$; 2) хоть одна точка пространства изолирована в бесконечном классе ее образов при всех преобразованиях группы (условие дискретности группы симметрии подобия).

Особая точка преобразования подобия называется особенной точкой группы симметрии подобия. Единственность особой точки в любой группе симметрии подобия не входит в определение 3. Она следует из теоремы, впервые доказанной Заморзаевым, которая формулируется следующим образом.

**Теорема** (Заморзаева [9]). Всякое преобразование $n$-мерной группы симметрии подобия сохраняет неподвижной любую её особенную точку.

Кроме того, в работах Шубникова и Заморзаева дана классификация всех групп симметрии подобия физического евклидового пространства. Как уже отмечалось, любое преобразование подобия $P$ с коэффициентом $k \neq 1$ может быть записано в виде произведения $P = K \cdot s$, $K$ – гомотетия с центром в точке $O$ и коэффициентом $k$, а $s$ – движение такое, что $s(O) = O$. Соответственно всем видам движений $s$ преобразований $P$ на плоскости имеют место следующие три вида (в обозначениях Шубникова):

1) операция $K$ – гомотетия;

2) операция $L$ – спиральное движение: эта операция слагается из последовательного выполнения операций поворота $v$ фигуры вокруг точки $O$ на некоторый угол $\varphi$ и гомотетии $K$ т.е. $L = K \cdot v$;

3) операция $M$ – гомотетическое отражение $M = K \cdot m$, где $m$ – отражение от прямой, проходящей через точку $O$.

Для трехмерного случая:

1) гомотетия $K$;

2) спиральное движение $L = K \cdot v$, где $v$ – вращение фигуры на некоторый угол $\varphi$, вокруг оси проходящей через точку $O$, с направляющим единичным вектором $\vec{l}$;

3) гомотетическое отражение $M = K \cdot m$, $m$ – отражение от плоскости, при этом $M^2 = K^2$, где $K^2$ имеет коэффициент подобия $k^2$;

4) спиральное отражение $\overline{L} = K \cdot \overline{v}$, где $\overline{v} = m \cdot v = v \cdot m$ – поворотное отражение, при этом $\overline{L} = L \cdot m$, $\overline{L}^2 = L^2$, причем коэффициент подобия $k^2$, ось $\vec{l}$ и угол $2\varphi$.

Отметим, что как топологические группы, группы симметрии подобия, удовлетворяющие определению 3, являются дискретными подгруппами непрерывных топологических групп, которые образуют обычные преобразования подобия, т.е. <u>группа симметрии подобия является подгруппой преобразований подобия</u> снабженная дискретной топологией.

Больше информации о симметрии подобия и геометрии Лобачевского можно найти, например, в книгах [10, 11].

## II. КОНФОРМНЫЕ ОТОБРАЖЕНИЯ

При геометрическом истолковании функции комплексного переменного $w(z)$ можно также рассматривать $w$ и $z$ как точки комплексной плоскости. Тогда функция $w(z)$ определяет некоторое соответствие между точками плоскости $z$ и плоскости $w$. При этом каждая геометрическая фигура в плоскости $z$ отображается в некоторую, вообще говоря, другую фигуру в плоскости $w$. При отображении, осуществляемом аналитической функцией сохраняются не только углы, но и направление обхода. Однако существуют отображения, которые сохраняют углы, но изменяют направление обхода. Примером может служить зеркальное отражение относительно оси $\mathrm{Re}(z)$, определяемое функцией $w = \bar{z}$. Такие отображения называются антиконформными (или конформными отображениями второго рода). Последовательное применение двух антиконформных отображений даёт конформное отображение. Каждая аналитическая функция определяет конформное отображение, при этом, обратно, каждое конформное отображение, т.е. отображение, сохраняющее углы и направление обхода, определяется некоторой аналитической функцией. Для того чтобы преобразование было конформным, необходимо, чтобы преобразование линейных элементов было преобразованием подобия. Отметим, что суперпозиция аналитических функций есть снова аналитическая функция.

Таким образом, преобразование плоскости, осуществляемое аналитической функцией, обладает следующим важным свойством в окрестности точки $z$. Бесконечно малые векторы всех направлений, выходящие из этой точки: 1) увеличиваются (или уменьшаются) по своей длине в одно и то же число раз (с точностью до бесконечно малых высшего порядка); 2) поворачиваются на один и тот же угол. Таким образом, фигуры в бесконечно малой области преобразуются в себе подобные – сохраняют форму. Фигуры конечных размеров искажаются, но углы между кривыми сохраняются (конформность или консерватизм углов). В частности, координатные линии $x = \text{const}$ и $y = \text{const}$ в конформном отображении преобразуются в два семейства взаимно ортогональных кривых. Поэтому, с помощью аналитических функций можно получить множество прямоугольных систем криволинейных координат.

Далее приведены некоторые конформные отображения, а именно:

$$w(z) = z^2 = Re(z)\cdot Im(z) + i\cdot[Re^2(z) - Im^2(z)]/2; \qquad (1)$$

$$w(z) = 1/z = [Re(z) - i\cdot Im(z)]/|z|^2. \qquad (2)$$

Заметим, что отображение (1) эквивалентно переходу из декартовой системы координат в параболическую систему координат (в которой координатные линии являются конфокальными параболами), а отображение (2) соответствует переходу в систему координат обратных радиусов (инверсия). Ньютоновское векторное поле центральных сил тяготения можно также представить на плоскости в виде $w(z) = \dfrac{1}{z} = \dfrac{x - i\,y}{x^2 + y^2}$. Отображение (2) равносильно последовательному выполнению: отражения от мнимой оси $z + \bar{z} = 0$ и инверсии относительно окружности с центром в начале координат и радиусом, равным единице $z\bar{z} = 1$. Напомним, что инверсия относительно окружности является инволюцией, конформным отображением второго рода (данная функция комплексного переменного является антиголоморфной, откуда следует конформность инверсии). Такое преобразование вводит бесконечность в теорию комплексных величин; точка $z = 0$ соответствует бесконечно удаленной точке $w = \infty$ плоскости $w$; обратно, точке $z = \infty$ соответствует точка $w = 0$. Отметим, что отображение (2) является частным случаем дробно-линейного преобразования (гомографического преобразования) вида $w(z) = \dfrac{az + b}{cz + d}$, где $a$, $b$, $c$, $d$ – постоянные целые числа, причем детерминант $a \cdot d - c \cdot b = 1$. Дробно-линейная функция – это единственная общая функция, описывающая однозначное отображение всей $z$-плоскости, включая бесконечно удаленные точки, на всю $w$-плоскость. Она представляется тремя парами значений $z$, $w$, переводит окружности в окружности (прямую можно рассматривать как окружность с центром в бесконечно удаленной точке); ее инвариантом является двойное отношение четырех величин (более детально про конформные отображения можно найти в [12]). Для трехмерного случая, по индукции можно покоординатно записать инверсию от сферы (причем не обязательно единичного радиуса).

Дробно-линейные отображения образовывают группу относительно операции суперпозиции (группа автоморфизмов сферы Римана, называется также группой Мёбиуса). Эта группа является также комплексно-трехмерной группой Ли. Формально группа Мёбиуса является проективизацией группы $GL_2(C)$, то есть имеет место эпиморфизм $\begin{pmatrix} a & b \\ c & d \end{pmatrix} \to \dfrac{az + b}{cz + d}$. Группа Мёбиуса изоморфна также специальной ортохронной группе Лоренца $SO^{\uparrow}(1, 3)$. Всякое дробно-линейное отображение может быть представлено в виде комбинации: сдвигов, инверсий, поворотов и растяжений. То есть, $w(z) = w_4(w_3(w_2(w_1(z))))$, где $w_1(z) = z + \dfrac{d}{c}$, $w_2(z) = \dfrac{1}{z}$, $w_3(z) = -\dfrac{ad - bc}{c^2} z$, $w_4(z) = z + \dfrac{a}{c}$. Как для сферы Римана, так и для единичного круга дробно-линейными функциями исчерпываются все конформные автоморфизмы.

Группы линейных преобразований можно рассматривать над произвольными полями – конечными или бесконечными, дискретными или непрерывными [13]. Все $n$-мерные решетки аффинно-эквивалентны. Другими словами, аффинное $n$-мерное пространство содержит по существу только одну решетку.

Интерпретируя комплексное число $z = x + iy$ как точку $(x, y)$ евклидовой (точнее, конформной) плоскости, заметим, что окружность с центром на оси $x$ описывается уравнением вида $A(x^2 + y^2) + 2Bx + C = 0$ или $Az\bar{z} + B(k + \bar{z}) + C = 0$ или $\left(z + \dfrac{B}{A}\right)\left(\bar{z} + \dfrac{B}{A}\right) = \dfrac{B^2 - AC}{A^2}$, где $A, B, C$ – действительные числа, $A > 0$. Её центр имеет координаты $\left(\dfrac{B}{A}, 0\right)$, а радиус равен $\sqrt{B^2 - AC}$. Инверсия относительно окружности меняет местами точки $z$ и $w$, где $\left(w + \dfrac{B}{A}\right)\left(\bar{z} + \dfrac{B}{A}\right) = \dfrac{B^2 - AC}{A^2}$, так что $1 \cdot w\bar{z} + B(w + \bar{z}) + C = 0$ и $w = -\dfrac{B\bar{z} + C}{A\bar{z} + B}$. Поскольку $B^2 - AC > 0$, можно так подправить коэффициенты уравнения, что будет выполняться равенство $B^2 - AC = 1$.

В предельном случае, когда $A = 0$ (и, скажем, $B = -1$) получается формула $w = -\bar{z} + C$ определяющая отражение относительно прямой $z + \bar{z} = 0$ или $x = C/2$. Всякое дробно-линейное преобразование $w = \dfrac{az + b}{cz + d}$ с действительными коэффициентами, можно выразить в виде произведения такой инверсии и отражения (или двух отражений). В самом деле

$$w = -\frac{d\bar{z} + (d^2 - 1)/c}{c\bar{z} + d}, \quad w' = -w + \frac{a-d}{c}, \quad (c \neq 0)$$ следует, что $w' = \frac{az + (ad-1)/c}{cz + d}$, а из соотношения $w = -\bar{z}$ и $w' = -\overline{w} + b$ следует, что $w' = z + b$. Преобразованиям $w = -\frac{B\bar{z} + C}{A\bar{z} + B}$ и $w = -\bar{z} + C$ соответствуют инволютивные матрицы $\begin{pmatrix} -B & A \\ -C & B \end{pmatrix}$, $\begin{pmatrix} -1 & 0 \\ C & 1 \end{pmatrix}$ с определителем –1. В частности, три преобразования $w = \frac{1}{\bar{z}}$, $w = -\bar{z} + 1$, $w = -\bar{z}$, соответствующие порождающим $R_1 = \begin{pmatrix} 0 & 1 \\ 1 & 0 \end{pmatrix}$, $R_2 = \begin{pmatrix} -1 & 0 \\ 1 & 1 \end{pmatrix}$, $R_3 = \begin{pmatrix} -1 & 0 \\ 0 & 1 \end{pmatrix}$ группы ß$_2$, которую они порождают.

Полуплоскость $y > 0$ можно рассматривать как конформную модель гиперболической плоскости. Прямые из этой плоскости изображаются полупрямыми и полуокружностями, ортогональными к оси $x$. Указанные выше преобразования можно интерпретировать в гиперболической плоскости как отражения относительно сторон треугольника с углами 0, π/2, π/3. Отсюда следует, что группа ß$_2$ – это группа [3, ∞], а её генетический код $R_1^2 = R_2^2 = R_3^2 = (R_1 R_2)^2 = (R_1 R_3)^2 = E$.

Основная задача настоящей работы показать, что для всех дискретных групп симметрии подобия существует разбиение плоскости, имеющее группы симметрии подобия каждого типа классификации Шубникова-Заморзаева. Из них особо выделим квазикристаллографические группы симметрии подобия, имеющие, в том числе поворотные симметрии 5, 8, 10 и 12-го порядков, которыми обладают квазикристаллы. При этом приведем некоторые исследования по конформным отображениям, которые позволяют расширить представления о подобии.

# III. АЛГОРИТМ ПОСТРОЕНИЯ КВАЗИРЕШЕТОК

Вначале матрично-векторным способом были выведены уравнения в рамках теории Новикова (подробнее см. работу [4]), позволяющие непосредственно рассчитывать координаты атомов в квазирешетках (кристаллических множеств обладающих поворотными симметриями 5, 6, 8, 10, 12, 14, 16, 18, … порядков) и изображать их на комплексной плоскости. Под словом «атом» подразумевается шар, центр, которого совпадает с элементом из кристаллического множества. Применяя целые алгебраические числа: $\rho = 1 + \sqrt{2}$; $\tau = (1 + \sqrt{5})/2$; $\eta = 1 + \sqrt{3}$, некоторые из уравнений (а именно для квазирешетки с ротационной симметрией: 8-го порядка, используя число $\rho$; 5-го и 10-го порядка – число $\tau$; 12-го порядка – число $\eta$), приведены в компактном виде в работе [4].

Целые алгебраические числа $\rho$, $\tau$, $\eta$ в литературе именуются металлическими сечениями: $\rho$ – число серебряного сечения, $\tau$ – число золотого сечения, $\eta$ – число бронзового сечения. Также они известны под названием чисел Пизо. Числа $\rho$, $\tau$, $\eta$ являются корнями квадратных уравнений $x^2 - 2x - 1 = 0$, $x^2 - x - 1 = 0$, $x^2 - 2x - 2 = 0$ соответственно. Отметим, что каждое из чисел $\rho$, $\tau$, $\eta$ является образующим элементом бесконечной циклической группы: $(\tau)_\infty = \{…, \tau^0 = 1, \tau^1 = \tau, \tau^2 = \tau + 1, \tau^3 = 2\tau + 1, \tau^4 = 3\tau + 2, \tau^5 = 5\tau + 3, \tau^6 = 8\tau + 5, …\}$; $(\rho)_\infty = \{…, \rho^0 = 1, \rho^1 = \rho, \rho^2 = 2\rho + 1, \rho^3 = 5\rho + 2, \rho^4 = 12\rho + 5, \rho^5 = 29\rho + 12, \rho^6 = 70\rho + 29, …\}$; $(\eta)_\infty = \{…, \eta^0 = 1, \eta^1 = \eta, \eta^2 = 2\eta + 2, \eta^3 = 6\eta + 4, \eta^4 = 16\eta + 12, \eta^5 = 44\eta + 32, \eta^6 = 120\eta + 88, …\}$. Заметим, что коэффициенты при $\tau$ и единице в выражениях для неотрицательных степеней $\tau$ в группе $(\tau)_\infty$ являются числами из последовательности Фибоначчи (0, 1, 1, 2, 3, 5, 8, 13, 21, 34, 55…). По аналогии запишем первые элементы из последовательности Пелля (0, 1, 2, 5, 12, 29, 70, 169, 408, 985, …) и последовательности (0, 1, 2, 6, 16, 44, 120, 328, 896, 2448, …), которые генерируются целыми алгебраическими числами ρ и η соответственно.

## IV. ВИЗУАЛИЗАЦИЯ И АНАЛИЗ ДИСКРЕТНЫХ ГРУПП СИММЕТРИИ ПОДОБИЯ ДЕЙСТВУЮЩИХ НА ПЛОСКОСТИ

Выполнение разбиений фигур, позволяющих визуализировать действие дискретных групп симметрии подобия на двумерных решетках и квазирешетках (напомним, что подгруппа $G$ группы $E_n$ всех движений $n$ – мерного евклидова пространства $E^n$ называется $n$ – мерной квазикристаллографической группой (в смысле Новикова), если ее пересечение с подгруппой $R^n \subset E_n$ всех трансляций, есть некоторая квазирешетка $T \subset R^n$), производилась следующим образом. Вначале строилось конкретное кристаллическое множество, обладающее поворотной симметрией 6, 8, 10, 12, … порядка. Далее из геометрического центра (особой точки) данного множества проводилось определенное число (равное порядку циклической группы) радиально расходящихся лучей. Затем каждый такой луч соединялся отрезком с соседним лучом под соответствующим углом и т.д. Такие углы рассчитывались исходя из конкретного расположения элементов множества и группы его симметрии. Напоследок осуществлялась различная цветовая раскраска полученных гомотетических тайлов.

Итак, исследуем, как выглядит объединение кристаллического множества, наложенного само на себя надлежащим образом до и после применения к нему гомотетии.

На рис. 1, а приведено объединение элементов кристаллического множества (точечная группа симметрии – диэдральная группа $D_6$) до и после осуществления операции $K$, с коэффициентом гомотетии $k = 2$. Исходное кристаллическое множество моделирует фрагмент двумерного гексагонального кристалла.

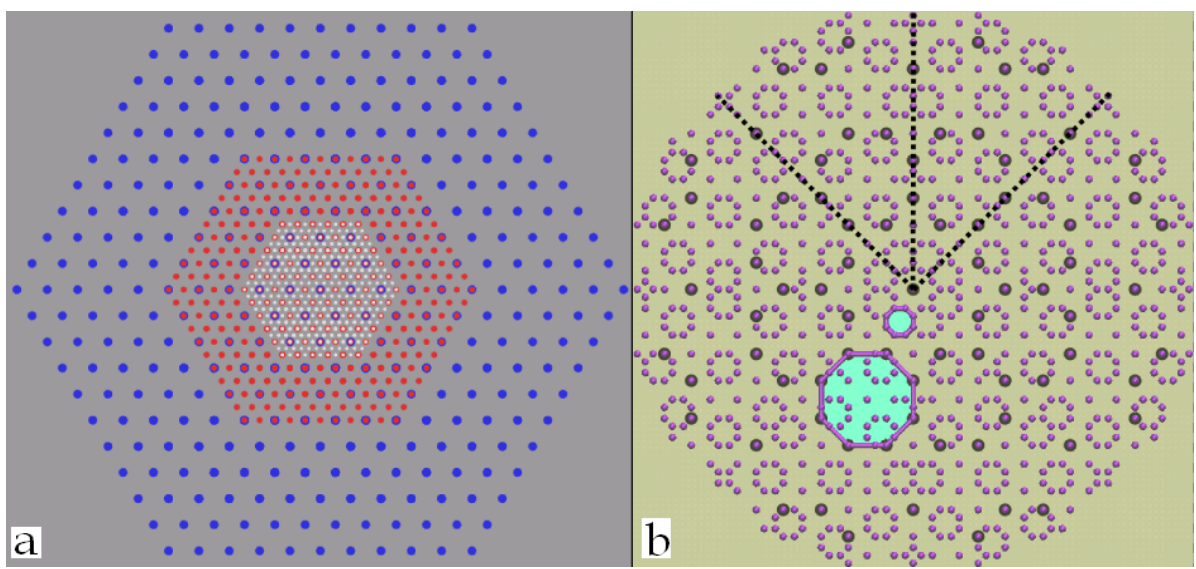

Рис. 1. Объединение кристаллического множеств до и после выполнения операций гомотетии: на a – результирующее объединение (как и кристаллическое множество – прообраз) обладает поворотной симметрией 6-го порядка; кружочки меньшего радиуса (белого цвета) соответствуют первоначальному расположению элементов кристаллического множества; кружочки среднего (красного цвета) и большего радиуса (синего цвета) отвечают положению элементов после одного и двойного осуществления операции $K$ соответственно; на b – результирующее объединение обладает поворотной симметрией 8-го порядка; кружочки меньшего радиуса соответствуют – первоначальному расположению элементов кристаллического множества; кружочки большего радиуса (изображена только часть этих элементов, т.е. только те, которые попали в окружность заданного радиуса) отвечают положению элементов после осуществления преобразования $K$.

На рис. 1, b изображено объединение элементов кристаллического множества (точечная группа симметрии – диэдральная группа $D_8$) до и после выполнения операции *K*, с коэффициентом гомотетии $k = 1 + \sqrt{2}$ (заметим, что $k = \rho$). Кристаллическое множество, к которому применяли гомотетию, моделирует двумерный октагональный квазикристалл.

На рис. 1, b пунктирными линиями выделены смежные «сектора». Элементы находящиеся как внутри, так и на границах выделенных секторов переходят друг в друга при поворотах на углы $\alpha = 2\pi n / 8$ (где *n* = 1, 2, …, 8) вокруг особой точки *O*, а также при отражениях от прямых и при инверсии в точке *O*. Минимальный угол поворота относительно точки *O*, при котором все множество переходит в себя, очевидно равен отношению 360° к порядку максимальной циклической подгруппы, соответствующей точечной группы.

Отметим, что каждый угол (выделенный пунктирными линиями) является фундаментальной областью циклической группы восьмого порядка. Также на рис. 1, b выделены два подобных октагона. Применяя операцию гомотетии к октагону меньшей площади получим октагон, который больше по площади к данному. Мы наблюдаем «процесс размножения» заданных частей фигуры в сторону их увеличения. Обратным построением можно получить части уменьшающиеся. Если принять во внимание, что кристаллическое множество – бесконечно, то после гомотетии «новые» положения элементов (образа) совпадают всюду со «старыми» (прообраза), другими словами – имеется масштабная инвариантность.

Таким образом, простейшая операция *K* является операцией переноса всех подобных частей фигуры в параллельное положение с одновременным увеличением (уменьшением) масштаба частей и расстояний между ними в *k* раз. При этом переносе соответственные точки подобных частей фигуры движутся прямолинейно.

На рис. 2, а изображена фигура, имеющая симметрию подобия 10 *L* ($\varphi = -\pi/5$), которая может совмещаться с собой как поворотами вокруг оси 10-го порядка на углы, кратные $2\pi/10$, так и спиральными поворотами на углы, кратные углу φ = –π/5. Знак минус перед значением угла показывает, что при спиральном движении от центра к периферии приходится совершать повороты по часовой стрелке. Поскольку группы рассматриваемого ряда содержат только операции первого рода (10 *L*), то фигура на рис. 2, a может существовать в двух энантиоморфных модификациях. Если половину частей данной фигуры раскрасить в черный цвет характерным образом (в том числе в шахматном порядке), то получаем новые фигуры с более низкой симметрией. На фигурах, приведенных, на рис. 2, следует отличать главные спирали от производных спиралей тем, что в них повторяющиеся части расположены наиболее густо.

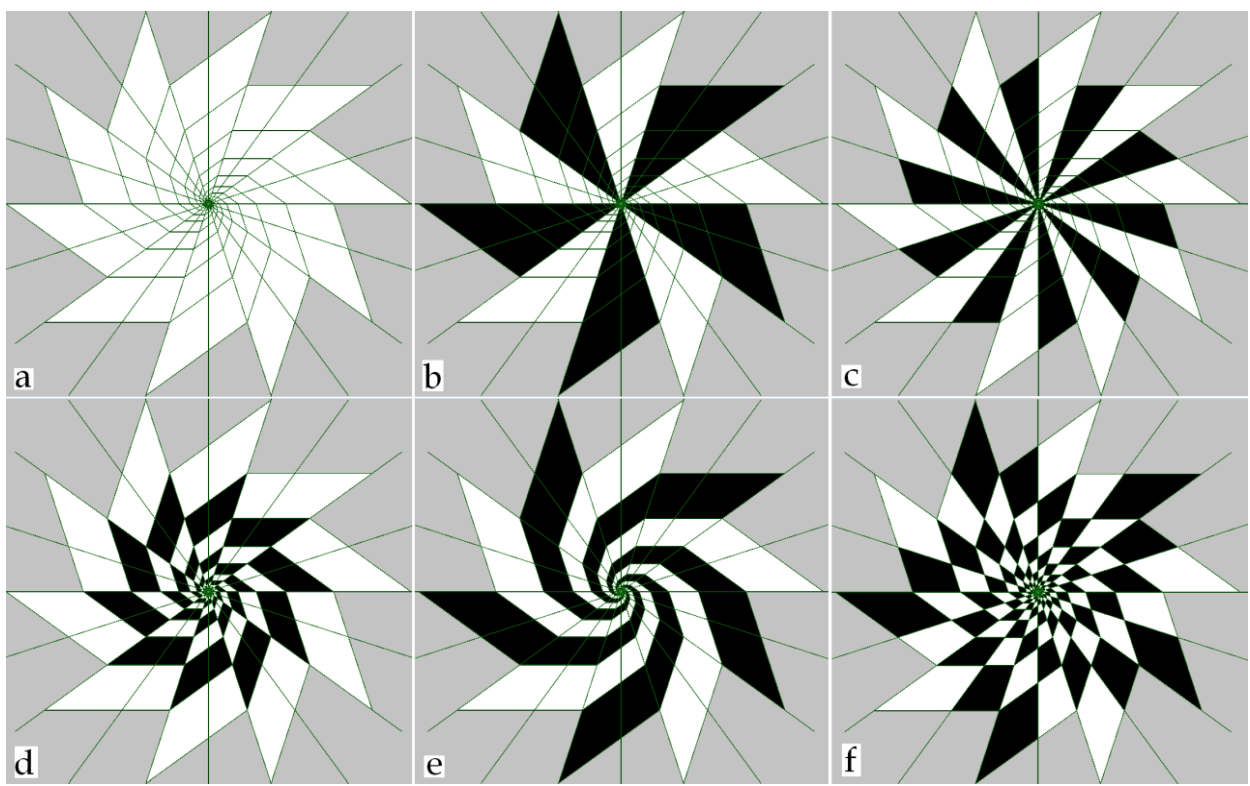

Рис. 2. Геометрические фигуры, обладающие группами симметрии подобия, подгруппами которых, являются циклические группы 5 и/или 10-го порядка.

На рис. 2, b изображена фигура, имеющая симметрию подобия 5 *L* (φ = –2π/5). На рис. 2, c изображена фигура, имеющая симметрию подобия

$10\,L$ ($\varphi = -\pi/5$). На рис. 2, d изображена фигура, имеющая симметрию подобия $10\,L$ ($\varphi = -\pi/5$). На рис. 2, e изображена фигура, имеющая симметрию подобия $5\,L$ ($\varphi = -\pi/5$).

На рис. 2, f изображена фигура совмещающаяся с собой при поворотах 5-го порядка.

В природе существуют физические объекты, геометрически равные, но взаимно полярные, т.е. которые имеют противоположный знак определенной физической характеристики. Такие объекты называют антисимметричными. Понятие антисимметрии введено Шубниковым в работе [14]. Преобразования, посредствам которых фигура переходит в себя, и при этом меняет цвет, называются операциями антисимметрии. Кроме того, в работе [15] было введено понятие цветной симметрии.

Группы симметрии подобия могут содержать комбинированные преобразования, например классическую группу и подгруппу антисимметрии некристаллографического (квазикристаллографического) типа (если преобразования перекрашивания эквивалентных частей фигуры рассматривать как операции антисимметрии).

На рис. 3 изображены различные варианты раскраски фигуры рис. 3, а, что позволяет визуализировать различные подгруппы её группы симметрии подобия.

На рис. 3, а дан пример фигуры, имеющей группу симметрии подобия $10\,m\,L$ ($\varphi = \pi/5$). Символ этой группы показывает, что фигура имеет одну ось симметрии десятого порядка, десять плоскостей симметрии и ось подобия с элементарным углом поворота, равным $1/10$ полного оборота.

Если на фигуре рис. 3, а раскрасить повторяющиеся части в шахматном порядке в два разных цвета, то получим фигуру, имеющую группу симметрии подобия $10\,m\,L$ ($\varphi = 2\pi/5$).

Изоэдральное разбиение фигуры на криволинейные гомотетические шестиугольники приведено на рис. 3, b. Данная фигура обладает группой симметрии 5 *m L* (φ = π/5).

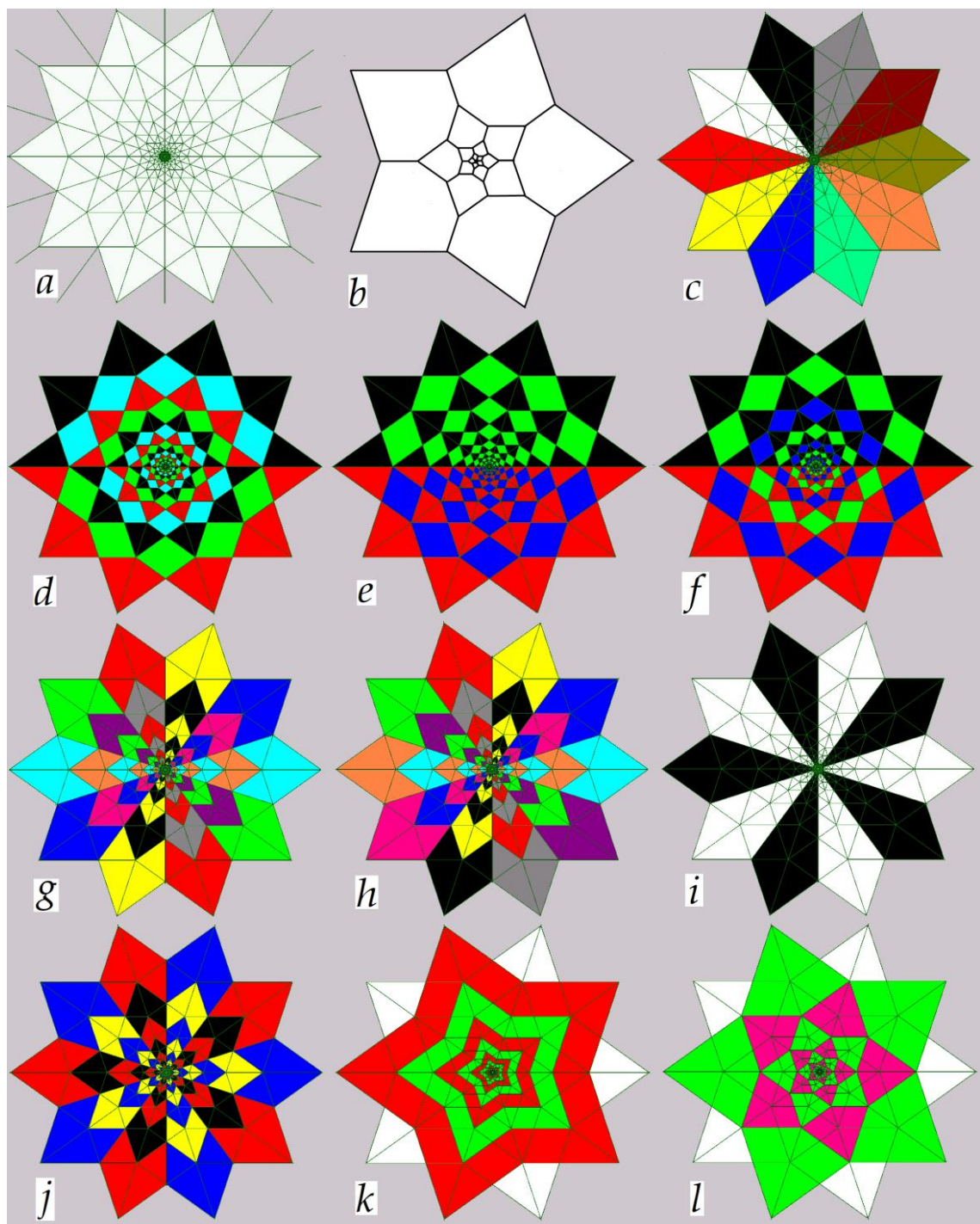

Рис. 3. Геометрические фигуры, обладающие группами симметрии подобия подгруппами которых являются циклические группы 1, 2, 5 и 10-го порядка; цветные фигуры приведенные на *b-l* получены из фигуры *a*.

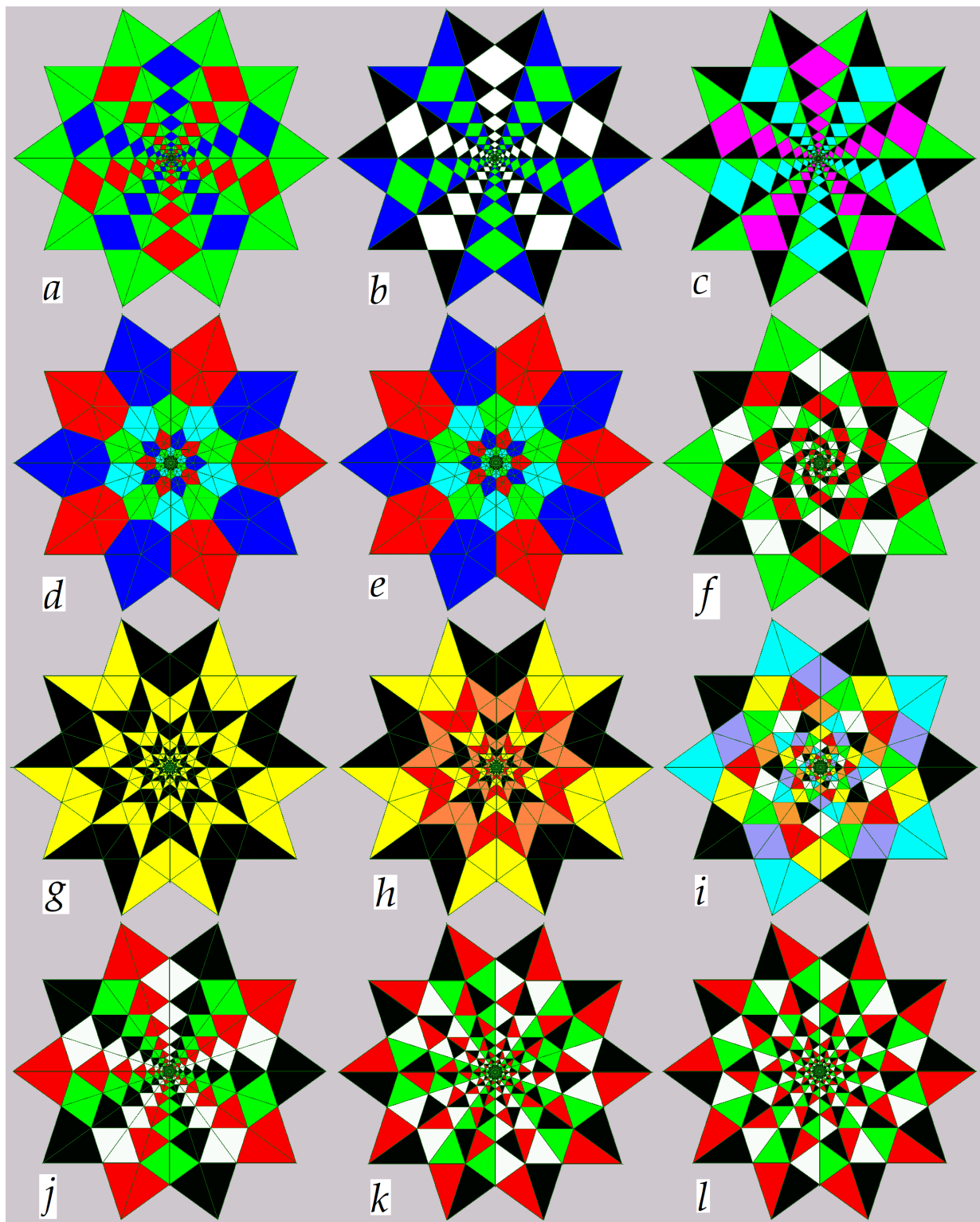

Рис. 4. Цветные фигуры приведенные на *a-l* получены из фигуры рис. 3, *a*; *f* – многозаходная спираль.

На рис. 4, *g* изображена фигура с группой симметрии 5·*mM*.

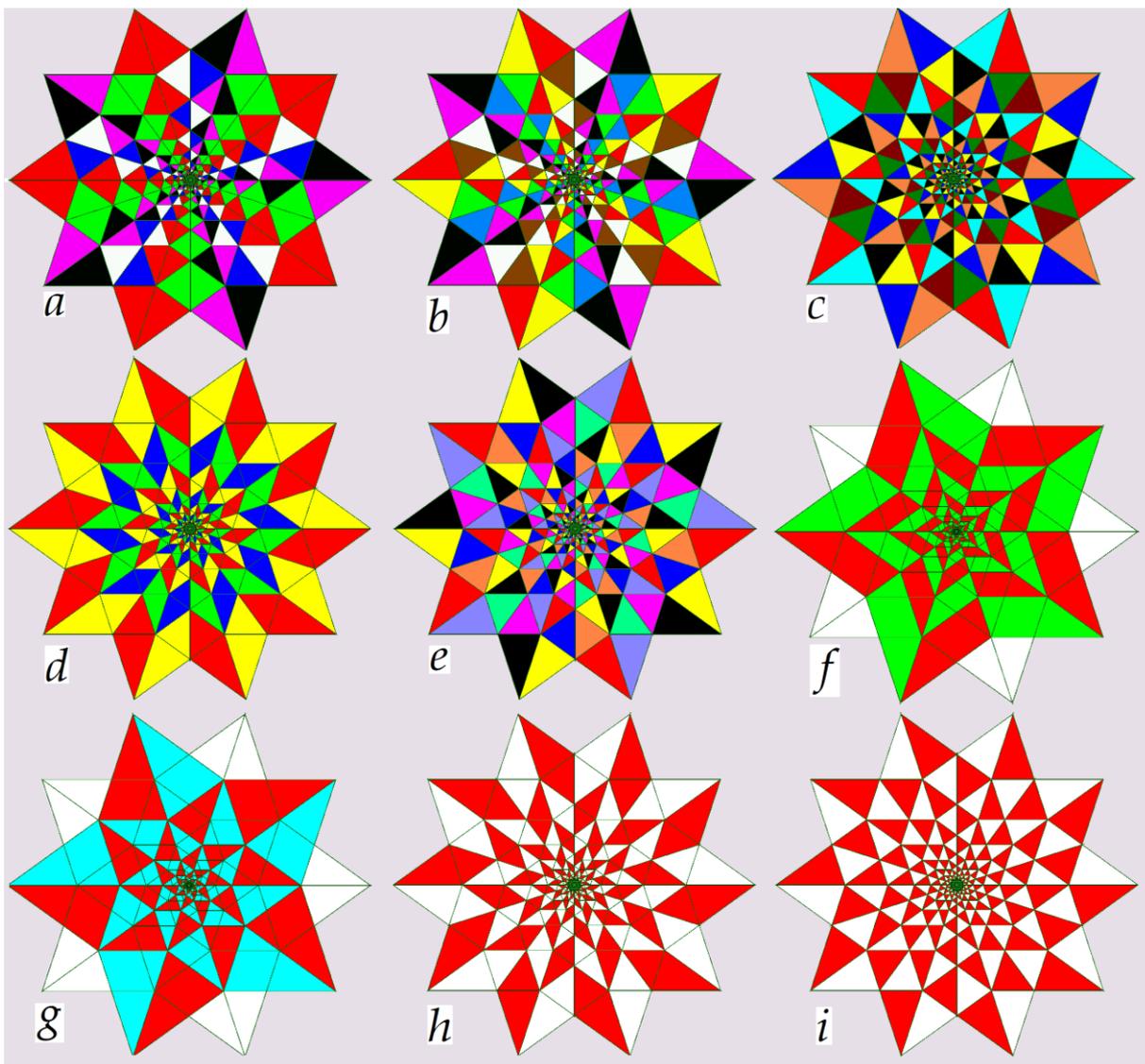

Рис. 5. фигуры приведенные на *a-i* получены из фигуры рис. 3, *a*.

На рис. 3, *f* изображена фигура с группой симметрии 5·*mM*.

Читателям предлагается в качестве упражнения определить группы симметрии подобия остальных двумерных фигур изображенных на рис. 3-6, а также трехмерных на рис. 0.

Отметим, что фигуры имеющие симметрии подобия, описываемые Шубниковым в работе [2], состоят из одного типа подобных плиток; фигуры же, приведенные на рис. 3-6, состоят из двух типов тайлов, что позволяет произвести, в том числе четырехцветную раскраску таких фигур.

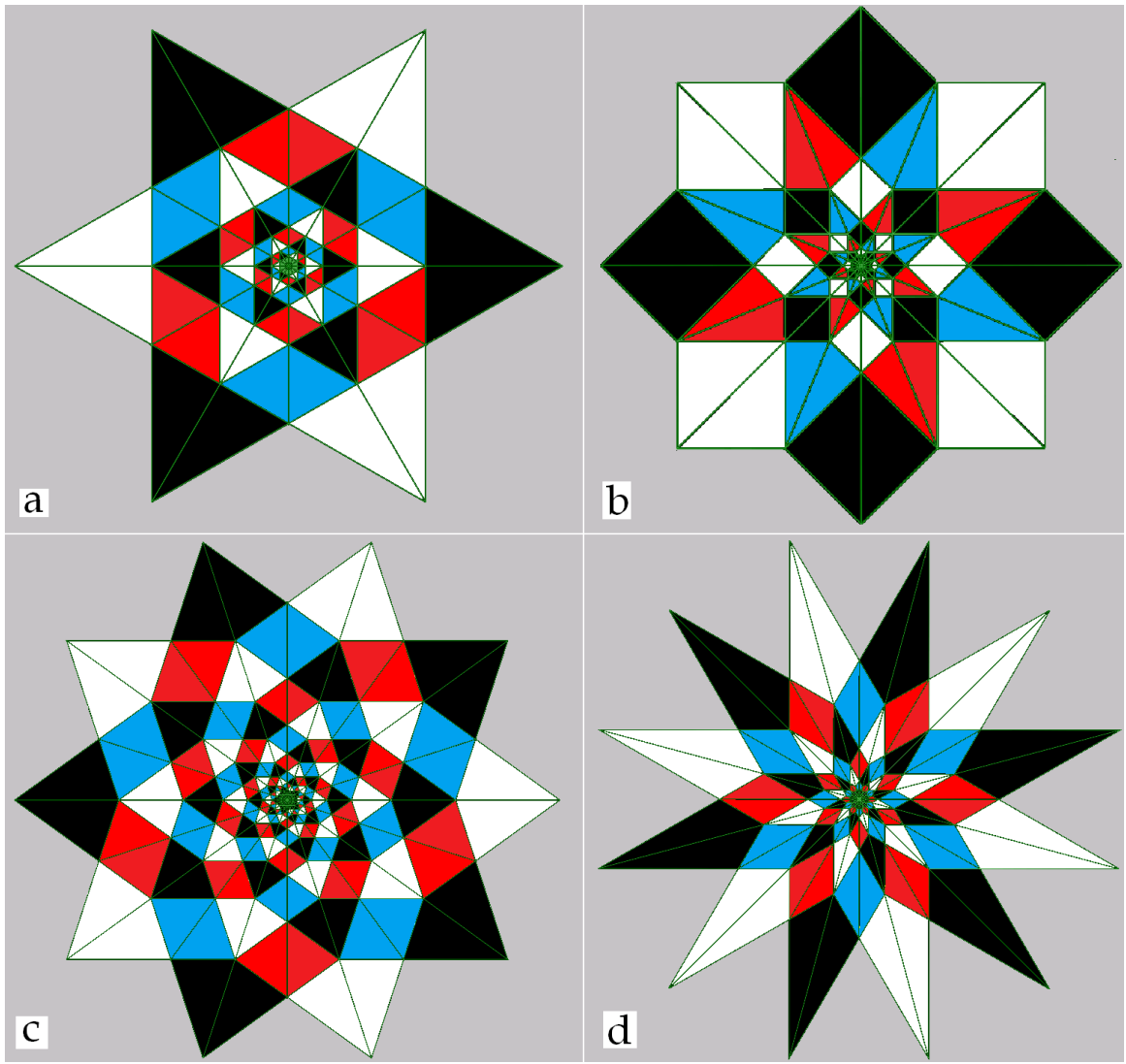

Рис. 6. Четырехцветные разбиения фигур, обладающих также поворотными осями симметрии 6, 8, 10 и 12-го порядка соответственно.

Некоторые из фигур, приведенных на рис. 3-6, получаются с помощью операции *K* и ее повторений, состоят из бесконечного множества частей, уменьшающихся при приближении к особой точке и увеличивающихся при удалении от нее. Уменьшающиеся части исчезают в особой точке; увеличивающиеся уходят в бесконечность. Прямые, проходящие через соответственные точки частей фигуры, сходятся в особой точке *O*. Причем у каждой фигуры число таких прямых равно числу равному порядку соответствующей диэдральной группы.

Далее перейдем к исследованию некоторых конформных отображений кристаллических множеств.

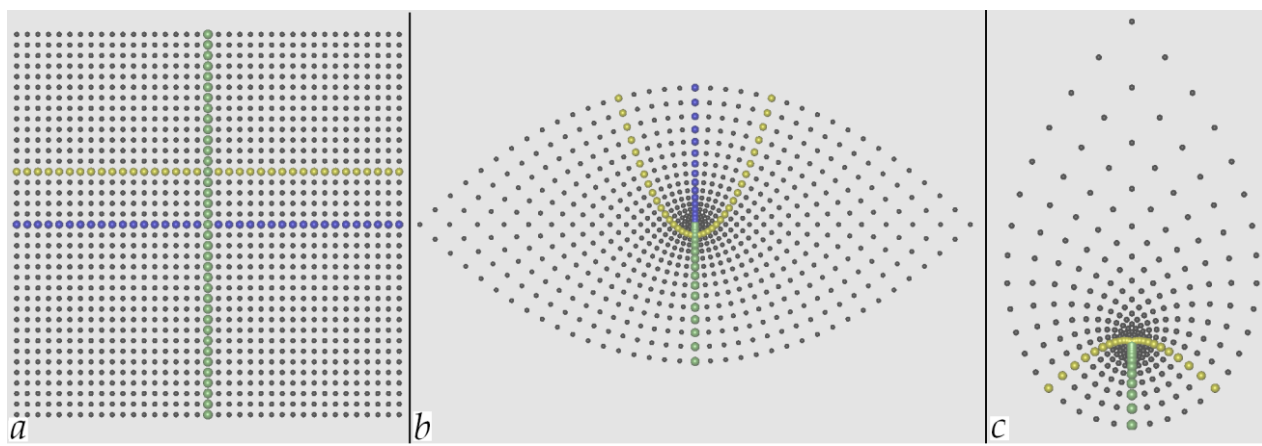

Рис. 7. *a* – двумерное кристаллическое множество, порядок циклической подгруппы точечной группы симметрий которого равен 4;
на *b* изображен образ кристаллического множества *a* после отображения (1);
*c* – образ кристаллического множества *a* после дважды примененного к нему преобразования (1).

На рис.7, *a* изображена модель квадратной решетки, фундаментальная область подгруппы трансляций которой является квадрат, а подгруппы поворотов – равнобедренный треугольник с углами 45°, 45°, 90°. Для визуализации особенностей при осуществлении отображения (1) «атомы» в некоторых рядах решетки изображены определенными размерами и окрашены характерными цветами. На рис. 7, *b* и *c* приведены образы множества (рис. 7, *a*), которые ограничены двуугольником и одноугольником соответственно. Угол при вершине у двуугольника и одноугольника равен 90°. Трехмерные аналоги множеств изображенных на рис. 7, *b* и *c*, согласно Скотту [16], являются орбифолды «чечевица» и «слеза» соответственно. Напомним, что *n*-мерным орбиобразием (орбифолдом) называется отделимое (хаусдорфово) паракомпактное пространство, локально гомеоморфное фактор-пространству $R^n$ по действию некоторой конечной группы.

Ряды атомов, проходящие через геометрический центр кристаллического множества (рис. 7, *a*) после применения к нему преобразования (1) переходят в полупрямые, выходящие из особой точки (рис. 7, *b* и *c*). После повторного осуществления преобразования (1) некоторые ряды «атомов» накладываются один на другой надлежащим образом. В центре масс множеств, приведенных рис. 7, *b* и *c*, расположена соответствующая особая точка (причем при $z = 0$ конформность нарушается).

Геометрическая модель гексагональной решетки, к которой было применено преобразование (1) и, следовательно, получена фигура (рис. 8, *a*), состояла из равносторонних треугольников (в вершинах которых располагались «атомы»), точечная группа симметрий которой – диэдральная группа $D_6$ (эта гексагональная решетка (прообраз) изображенная на рис. 1, *a*). Фигура, ограничивающая множество, изображенное на рис. 8, *a* – «треугольник», сумма внутренних углов у которого равна 360°. Отметим, что прямыми у фигур, изображенных на рис. 7 и рис. 8 являются конфокальные параболы, фокус которых совпадает с особой точкой. На рис. 8, *b* показан вариант двухцветной раскраски фигуры (рис. 8, *a*) в шахматном порядке, которая обладает поворотной симметрией 1-го порядка. На рис. 8, *c* приведен вариант «секторальной» двухцветной заливки фигуры (рис. 8, *a*) в шахматном порядке, которая обладает поворотной симметрией 3-го порядка. На рис. 8 *d-f* приведены еще некоторые варианты раскраски фигуры (рис. 8, *a*), обладающие поворотной симметрией 3-го порядка. Заметим, что количество подобных тайлов в разбиении (например, рис. 8, *d*) увеличивается при отдалении от особой точки. Из фигуры (рис. 8, *d*) следует, что центры подобных тайлов разбиения, шестиугольников лежат на чередующихся через одну конфокальных параболах, причем, чем ближе расположена вершина параболы к ее фокусу (т.е. чем больше кривизна параболы), тем больше на ней расположено центров шестиугольников.

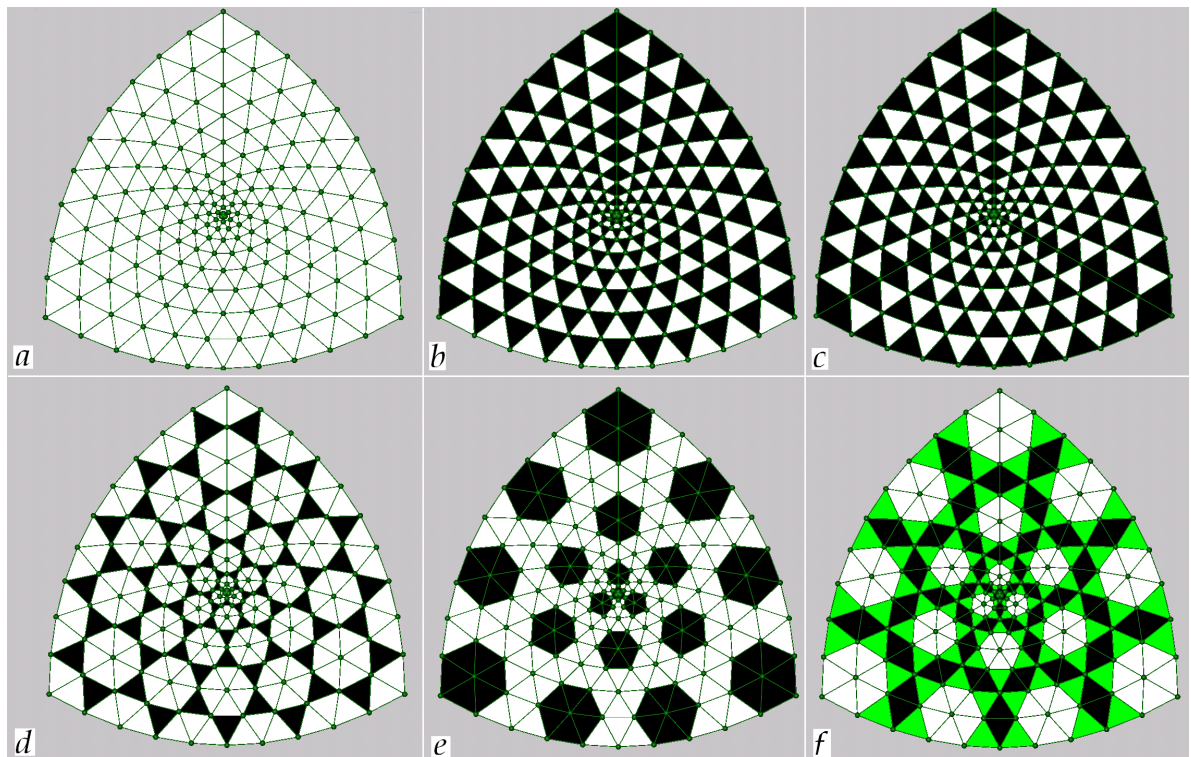

Рис. 8. *a* – образ двумерной гексагональной решетки, к которой применили преобразование (1) и соединили ближайшие элементы множества отрезками; на *b-f* приведены различные варианты цветовой раскраски фигуры *a*.

Фигуры изображенные на рис. 8, *a, c, d, e* и *d* обладают в том числе поворотной симметрией 3 порядка. Фигура на рис. 8, *b* совмещается с собой при тождественном повороте.

Рассмотрим инверсию от окружности, т. е. преобразование действительной оси посредством обратных радиусов, на примере квадратной решетки изображенной на рис. 7, *a*.

Прямыми на рис. 9 являются кусочно-гладкие полуокружности. В каждой вершине разбиения фигуры рис. 9, *a* сходятся по четыре «подобных» четырехугольника (за исключением границ фигуры). Четыре периферийных ряда элементов множества, приведенного на рис. 7, *a*, после применения к нему отображения (2) преобразуются в четыре полуокружности, разграничивающие особую область с разбиением фигуры рис. 9, *a*.

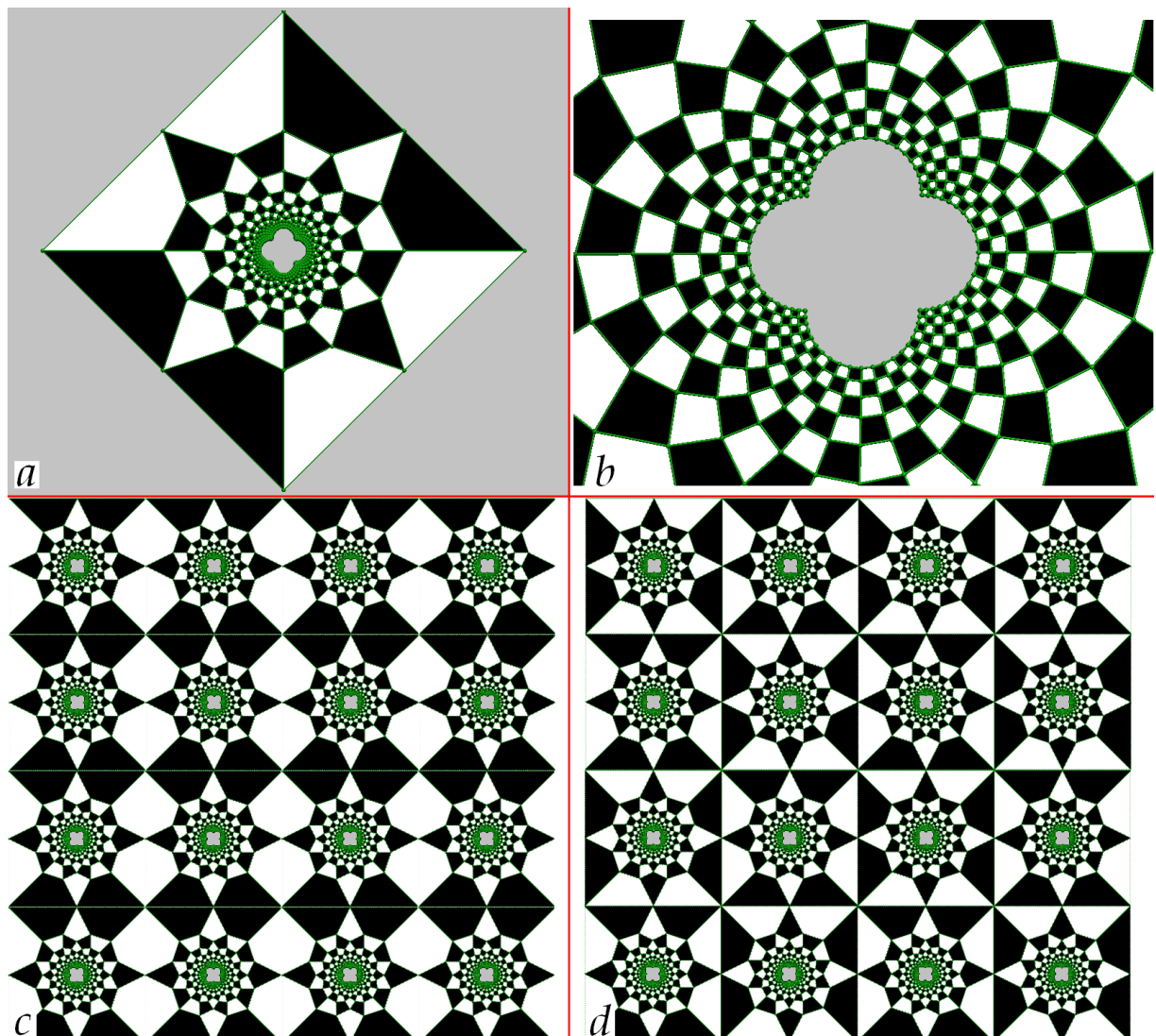

Рис. 9. *a* – фигура, полученная путем применения преобразования (2) к кристаллическому множеству, изображенному на рис. 7, *a* (каждый элемент соединен с ближайшим отрезком прямой); на *b* изображен увеличенный фрагмент центральной части фигуры *a*; на *c*, *d* показаны варианты разбиения плоскости грань в грань фигурой *a*, которые обладают поворотной симметрией 2-го и 4-го порядка соответственно.

Под особой областью подразумевается центральная часть фрагмента фигуры, изображенной на рис. 9, *b*. Четыре полуокружности,

ограничивающие особую область, пересекаются таким образом, что образуют четыре особые точки. Следовательно, можно заключить, что число подобного рода особых точек в образе любого кристаллического множества, которое может быть получено при осуществлении инверсии относительно окружности последнего, равно числу вершин многоугольника, ограничивающего данное множество. По индукции на примере кубической решетки в трехмерном пространстве после осуществления инверсии относительно сферы у образа будет 8 особых точек. Все 8 особых точек образуются путем пересечения 6 полусфер, причем каждая особая точка образуется пересечением 3 полусфер (см. рис. 10).

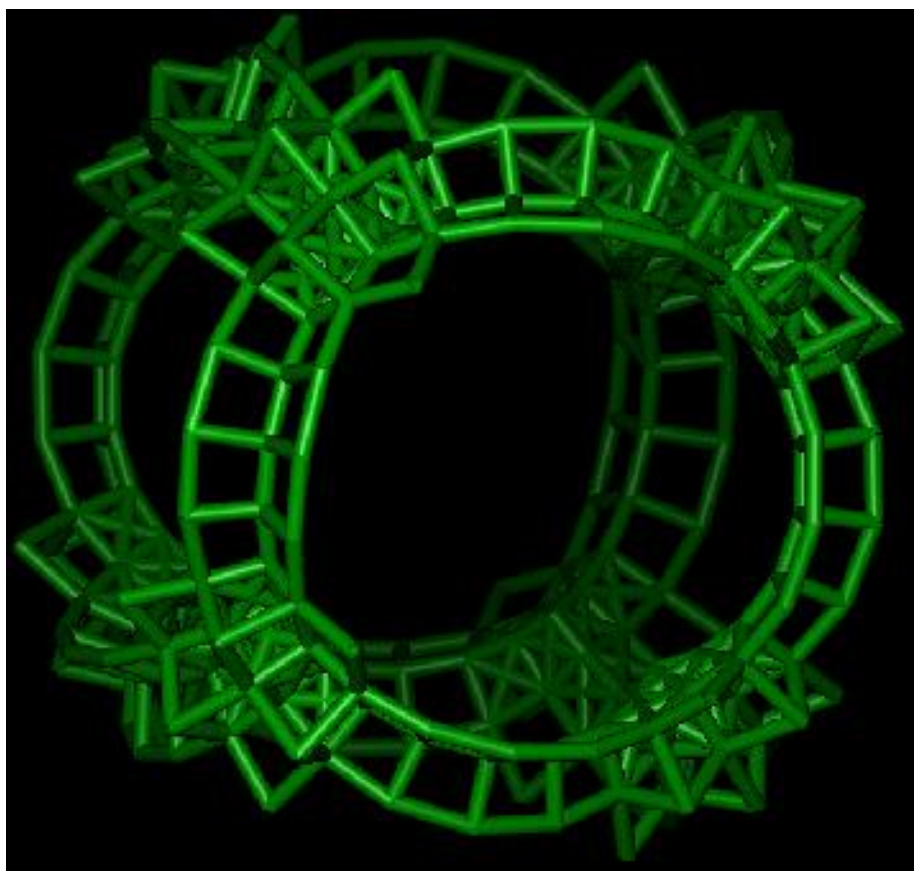

Рис.10. Фрагмент центральной части кристалла

(кубическая сингония *P*) в системе координат обратных радиусов;

элементы кристаллического множества для лучшей визуализации соединены отрезками прямых.

Так как фигура, изображенная на рис. 9, *a,* ограничена квадратом, то возможно произвести соответствующее разбиение плоскости «грань в грань». На рис. 9, *c* приведено разбиение плоскости фигурой рис. 9, *a* (повернутой на угол 45° против часовой стрелки), подействовавши на которую группой параллельных переносов. Разбиение плоскости, приведенное на рис. 9, *d*, получается путем последовательного выполнения трансляции фигуры рис. 9, *a* (вначале повернутой на угол 45° против часовой стрелки) и отражения ее в антизеркале (причем трансляция и отражение коммутируют).

Заметим, что количество «подобных» тайлов в «секторальном» разбиении (см. рис. 9 *b*) увеличивается при приближении к особой области.

Представляется весьма интересным сыграть в шахматы (шашки) на доске аналогичной фигуре изображенной на рис. 9, *a.*

Рассмотрим один важный вопрос, связанный с задачей рентгеноструктурного анализа кристаллов. Если возьмем все отражения от граней кристалла рентгеновских лучей и получим их стереографическую проекцию на плоскость, параллельную фотографической пластинке и проходящую через кристалл, то такая проекция называется стереоциклической. По сути дела она является не чем иным, как стереографической проекцией отраженных лучей, но для ребер кристалла она уже не будет обычной стереографической проекцией. Ребра кристалла изображаются при помощи сети кругов, проходящих через центр проекции. На таких диаграммах как бы складываются две проекции: стереографическая проекция отраженных лучей (пятен); циклическая (круговая) проекция ребер кристалла. Связь между Лауэграммой и стереоциклической проекцией такова же, как и связь между обычной гномонической и стереографической проекциями. Пятна, лежащие на каждом из кругов, являются стереографическими проекциями отраженных лучей от граней, принадлежащих одной и той же зане (поясу). Стереографическая проекция оси зоны не совпадает с центром соответствующей окружности, так как она является ее полюсом и отстоит от всех точек этой окружности на равную угловую величину. Малый круг, на котором лежат пятна, можно рассматривать как свободную проекцию оси (кругового конуса) некоторой зоны. Каждая такая окружность является геометрическим образом взаимно-однозначно связанным с осью пояса, которому принадлежат отражающие грани. Из этого следует, что экспериментальную стереоциклическую проекцию кристалла можно смоделировать путем построения соответствующей дуальной решетки и отображения ее в систему координат обратных радиусов с последующей ортогональной проекцией на плоскость.

Таким образом, инверсией комплексной плоскости $C$ относительно единичной окружности $П$ называется такое преобразование точек плоскости $C$, что точка $z$, находящаяся внутри окружности $П$, преобразуется в точку $w$, находящуюся на продолжении отрезка $Oz$ (где начало координат является центром инверсии) и расположенной вне окружности. причем произведение расстояний (в евклидовом смысле) от точки $O$ до отображенной и первоначальной точки равно единице. Отметим также определение инверсии, предложенное Тёрстоном в [17]. Пусть $П$ – окружность на евклидовой плоскости с центром в точке $O$. Инверсией $I$ относительно $П$ называется единственное отображение проколотой (в точке $O$) плоскости в себя, которое: оставляет все точки окружности $П$ неподвижными; меняет местами внутренность и внешность $П$; любую окружность, ортогональную $П$, переводит в себя. Это необходимое и достаточное условие существования инверсии. Точки $z$ и $w$, переходящие одна в другую с помощью инверсии, называются симметричными относительно окружности $П$. Аналитическое отображение $w(z)=1/z$, при $z \neq 0$ является движением первого рода. Данное отображение сохраняет углы во всех точках комплексной плоскости, в том числе и в случае $z=0$ и $z=\infty$. При этом углом двух линий при $z=\infty$ называется угол, образованный отображенными линиями посредством функции $1/z$ и плоскости $w$ при $w = 0$. Отметим, что инверсия $1/z$ является движением геометрии Лобачевского $H^2$, только при условии $y > 0$. Действительно $\dfrac{1}{\bar{z}} = \dfrac{z}{\bar{z} \cdot z} = \dfrac{x + iy}{x^2 + y^2}$, если $z = x + iy$ (мы избавились от иррациональности в знаменателе). Определитель матрицы этого преобразования равен $y(x^2 + y^2)$. Он больше нуля при $y > 0$, т.е. в верхней полуплоскости евклидовой плоскости.

Далее перейдем к рассмотрению конформных отображений кристаллических множеств моделирующих двумерные квазикристаллы.

На рис.11, *a*, *c* и *e* приведены множества элементов, моделирующие двумерные квазикристаллы, порядок циклических групп симметрий которых равен 14, 16 и 18 соответственно. В геометрическом центре каждой такой квазирешетки находится особая точка. Из образа кристаллического множества, например рис. 11, *b* следует, что средняя плотность элементов множества увеличивается при приближении к особой точке. Форма тайлов разбиения (тайлы получаются путем соединения отрезками каждого элемента множества с его ближайшим окружением) претерпевает некоторую деформацию, однако число различных тайлов в разбиении (образа и прообраза) остается неизменным. Отметим что, одним и тем же конкретным набором тайлов можно произвести, как замощение (разбиение) плоскости, обладающее некоторой определенной группой симметрии, так и соответствующие замощения плоскости, обладающие симметриями подгрупп данной группы.

Значение функции Эйлера для квазирешеток изображенных на рис. 11. *a, c* и *e* соответственно равны φ(14) = 6, φ(16) = 8 и φ(18) = 6.

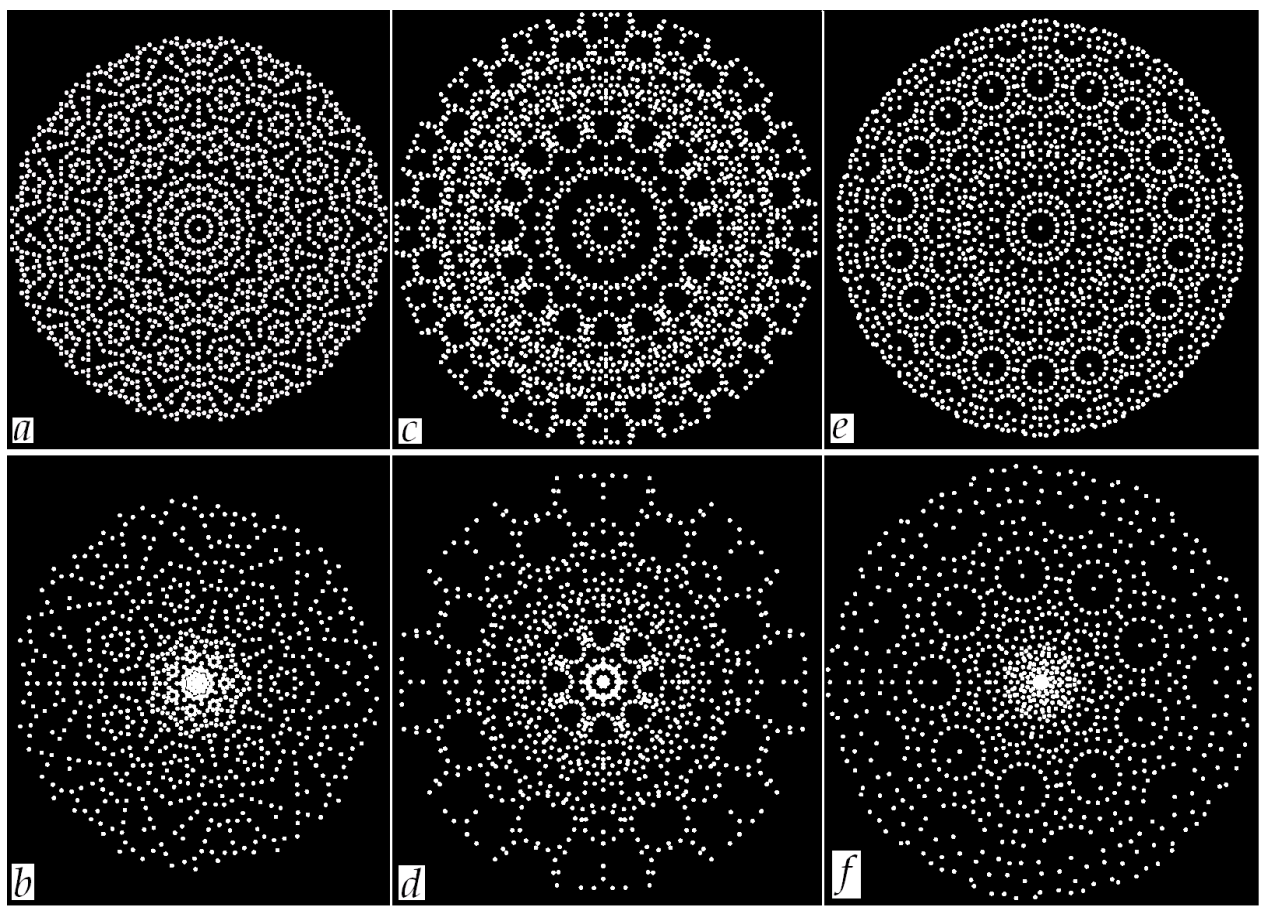

Рис. 11. На *a*, *c* и *e* изображены двумерные кристаллические множества (т.е. квазирешетки в смысле Новикова [18]), точечные группы симметрий, которых являются диэдральными группами $D_{14}$, $D_{16}$ и $D_{18}$ соответственно; *b*, *d* и *f* – образы множеств *a*, *c* и *e* соответственно, полученные после применения к последним преобразования (1).

Заметим, что когда порядок циклической подгруппы точечной группы симметрий кристаллического множества кратен двум, то после применения к этому множеству преобразования (1) порядок циклической подгруппы образа становится вдвое меньше по отношению к порядку циклической подгруппы прообраза. Такое уменьшения порядка циклической группы можно проследить на рис. 7 и рис. 12.

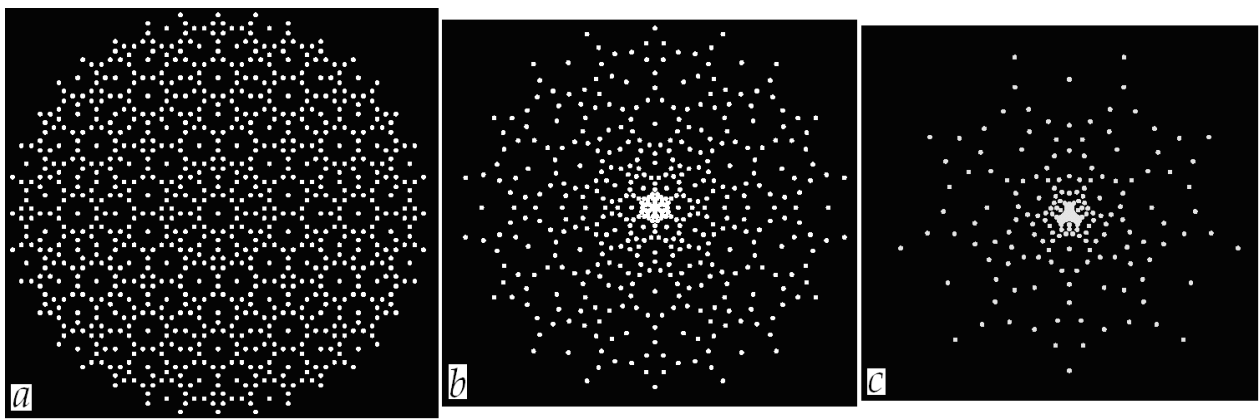

Рис. 12. На *b* и *c* приведены образы множества *a*, к которому последовательно применении конформное отображения (2) единожды и дважды соответственно.

Шубников в [19] пишет, что намеченная еще в работах Виола (1904) и Вульфа (1909) идея вывода групп аффинных деформаций получила дальнейшее развитие в работах Михеева (1961), Навилкина (1951), Дубова (1970), и Заболотного (1973) в форме так называемых групп гомологии и групп криволинейной симметрии. Хотя вывод этих групп, изоморфных цветным группам, не может считаться еще до конца завершенным, идея аффинных деформаций, восходящая к теории кристаллографических пределов Е. С. Федорова, открывает путь изучения динамической симметрии кристаллов и трактовки классических групп как усредненных по времени динамических групп. Группы по модулю локально-аффинных преобразований могут быть использованы при анализе симметрии геометрически неоднородных объектов.

По самому общему определению, данному Г. Вейлем, теория симметрии совпадает с теорией групп автоморфных преобразований. В ее основе лежит аксиома равенства. По сути дела отсюда и вытекает свойство преобразований, сохраняющих некий объект инвариантным, образовывать группу. Отсюда же следует и широта приложений теории симметрии, и, равным образом, очерчиваются ее границы, поскольку преобразования автоморфизма не исчерпывают всех отношений эквивалентности и порядка между объектами.

Любые обобщения теории симметрии сохраняют ее основные групповые постулаты, они идут по пути замены понятия равенства более широким понятием относительного равенства конкретных реализаций абстрактных групп, замены известных групп их изоморфными или гомоморфными представлениями или расширениями, по пути поисков новых автоморфизмов. Характеризуется научное значение таких поисков, достаточно сказать, что само обнаружение новых отношений эквивалентности и новых автоморфизмов у известных объектов уже означает проникновение на более глубокий структурный уровень исследования.

Стереографическая проекция не только конформна, но и обладает круговым свойством, заключающимся в том, что окружности и прямые плоскости отображаются в окружности на сфере, и обратно. Окружность на сфере можно получить как пересечение сферы с плоскостью. В отличие от проекции Меркатора стереографическая проекция дает взаимно однозначное отображение всей плоскости на поверхность сферы. Исключение составляет лишь южный полюс, которому не соответствует никакая точка плоскости. Это исключение можно устранить, присоединяя к собственным точкам плоскости бесконечно удаленную точку и считая, что она соответствует точке на сфере. Эта несобственная точка не изображается никаким комплексным числом, так как невозможно расширить семейство комплексных чисел введением «несобственного числа», не придя при этом в противоречие с правилами действий (например, выражения $0 \cdot \infty$, $\infty/\infty$ не могут быть однозначно определены). Однако на сфере можно рассматривать преобразования, являющиеся конформными и на южном полюсе, и, обратно, переводящие произвольную точку сферы в южный полюс и конформные в этой точке. Отметим, что рассмотрение преобразование сферы привело нас к представлению бесконечности на комплексной плоскости в виде одной точки, в отличие, от проективной геометрии на плоскости, в которой бесконечно удаленные точки образуют прямую линию. Такое представление бесконечности содержит некоторый элемент произвола. Отметим одно преимущество, что плоскость, как и сфера, становится замкнутой поверхностью и приобретает весьма простую топологическую структуру.

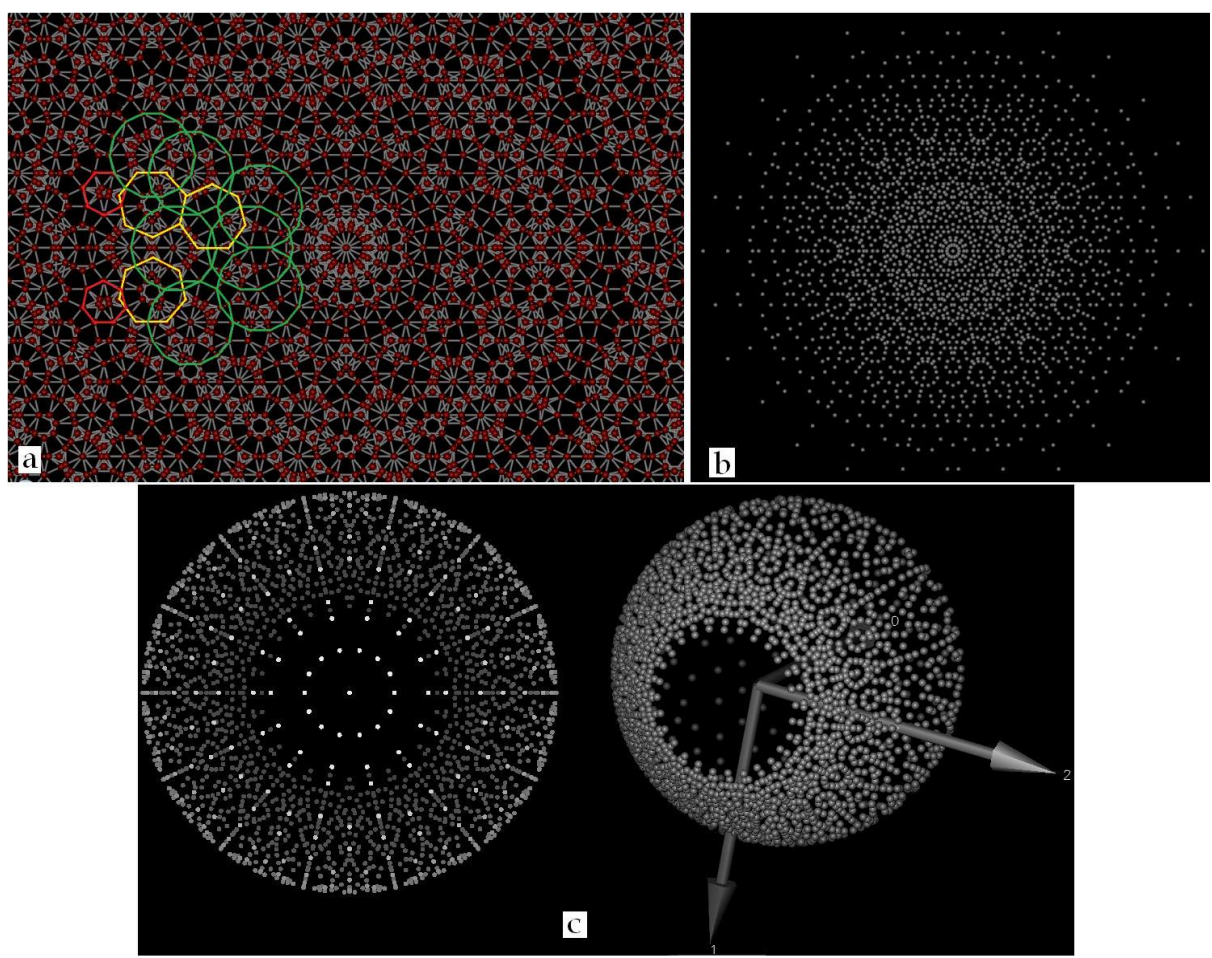

Рис. 13. Кристаллические множества с точечной группой симметрии $D_{14}$;

(c) – проекция на сферу Римана (в двух ракурсах)

кристаллического множества (b)

(рисунок взят из работы [4]).

На рис. 13, а приведено кристаллическое множество, обладающее ротационной симметрией 14-го порядка, на котором отдельно отрезками выделены четырнадцатиугольники и подобные гептагоны.

Преимущество построения квазирешеток, используя комплексные числа, заключается еще и в том, что можно осуществлять различные проективные преобразования, одним из которых является отображение на сферу Римана $\hat{C} = C \cup \{\infty\}$. Для соответствующего отображения применялись формулы стереографической проекции $x_0 = \frac{1-|z|^2}{1+|z|^2}$, $x_1 = \frac{2Re(z)}{1+|z|^2}$, $x_2 = \frac{2Im(z)}{1+|z|^2}$. Проективная сфера Римана имеет две особые точки: северному полюсу соответствует бесконечно удаленная точка, а южному полюсу – особая точка. Так как кристаллическое множество рис.13, b, спроецированное на сферу Римана рис.13, c, конечно, то элемент в северном полюсе отсутствует. По аналогии с квазипериодическими разбиениями плоскости произведены соответствующие разбиения сферы Римана. Конформное отображение комплексных чисел на сферу Римана применяется для построения квазикристаллов в моделях Пуанкаре пространства Лобачевского $H^2$. Применяя теорию вычетов, разработаны алгоритмы для разложения любого кристаллического множества с группой симметрии $D_m$ (при четном $m$) на подмножества, обладающие, в том числе поворотами порядка $m/2$.

## V. ЗАКЛЮЧЕНИЕ

Отправляясь от классических результатов Шубникова и Заморзаева, разработаны алгоритмы построения компьютерных моделей кристаллических множеств, которые обладают двумерными группами симметрия подобия классифицированными в работах данных исследователей.

Рассмотрены некоторые конформные отображения кристаллических множеств, которые имеют широкое применение в задачах математической физики, таких как: теория эффекта Штарка; теория кристаллического поля; задача о потенциале вблизи угла; электротехника; гидро- и аэродинамика.

Показано, что у «классических» двумерных кристаллов коэффициенты гомотетии являются целыми числами, а у квазикристаллов коэффициенты гомотетии – целые алгебраические числа квадратичного расширения поля рациональных чисел *Q*. Например, у двумерной решетки, обладающей поворотной симметрией 4 порядка (квадратная сетка) коэффициент гомотетии *k* равен 2, а у квазирешеток обладающих поворотными симметриями 5 и 10 порядка *k* = τ, для квазирешеток обладающих поворотными симметриями 8 и 12 порядка коэффициент *k* равен $1+\sqrt{2}$ и $1+\sqrt{3}$ соответственно.

Приведены разбиения фигур, которым присущи цветные группы симметрии подобия.

Установлено, что конформное отображение решетки кристалла в систему обратных радиусов эквивалентно построению стереоциклической проекции дуального кристалла к данному.

Если взять реальный кристалл, у которого одной из операций симметрии будет инверсионная ось 4-го (или 6-го) порядка (а именно, кристаллы с точечной группой симметрии: $\bar{4}$, $\bar{4}2m$, $\bar{4}3m$, $\bar{6}$ или $\bar{6}m2$) и получить от него экспериментальную лауэграмму снятой вдоль такой инверсионной оси, то в результате множество рефлексов на фотопластинке будет обладать поворотной симметрией восьмого (или двенадцатого) порядка! Этот факт можно объяснить следующим образом: лауэграмма – это проекция дуальной трехмерной решетки (в обратном пространстве) на сферу Эвальда; точечная группа симметрии кристалла и дуального к нему кристалла – изоморфны; так как лауэграмма это проекция на плоскость (сферу бесконечного радиуса), то очевидно, что теряется информация, о том была ли ось симметрии поворотной или инверсионной. Таким образом, если получить на лауэграмме множество рефлексов, которое будет обладать поворотной симметрией восьмого порядка, то это еще не означает, что объект, от которого она снята, является октагональным квазикристаллом!

Для того чтобы утверждать, что это октагональный квазикристалл необходимы дополнительные структурные исследования.

Следствием исследования можно считать, что спиральная тенденция в природе прослеживается и на нано уровне.

Результаты исследования могут иметь применение в декоративном искусстве. При этом имеют некоторые корреляции с символизмом в сакральной геометрии.

## VI. СПИСОК ЛИТЕРАТУРЫ